\documentclass[twoside]{article}
\usepackage[a4paper]{geometry}
\usepackage[latin1]{inputenc} % ou \usepackage[utf8]{inputenc}

\usepackage{tensor}
\usepackage{color}

\usepackage{graphicx}
\usepackage{ upgreek }
\usepackage{graphicx, amsmath,  amssymb}%, algorithmic}
\usepackage{accents}

\usepackage{amsthm}

\usepackage[ruled]{algorithm2e}
\usepackage{algorithmic}

\begin{document}

\newcommand{\revtex}{REV\TeX\ }
\newcommand{\classoption}[1]{\texttt{#1}}
\newcommand{\macro}[1]{\texttt{\textbackslash#1}}
\newcommand{\m}[1]{\macro{#1}}
\newcommand{\env}[1]{\texttt{#1}}

\newcommand{\C}{\mathcal{C}}

\newcommand{\SUIO}{$\mathcal{SUIO}$}
\newcommand{\B}{$\mathcal{B}$}
\newcommand{\VS}{$\mathcal{V}$}
\newcommand{\IMU}{$\mathcal{I}$}
\newcommand{\tIMU}{\textnormal{\IMU}}
\newcommand{\tVS}{\textnormal{\VS}}
\newcommand{\tih}{\widetilde{h}}
\newcommand{\tobs}{\widetilde{\OBS}}
\newcommand{\NO}{s}
\newcommand{\ND}{r}

\newcommand{\M}{\mathcal{M}}
\newcommand{\Li}{\mathcal{L}}
\newcommand{\RM}{\mathcal{RM}}
\newcommand{\Hf}{\mathcal{H}}

\newcommand{\oplusn}{+}
\newcommand{\bigoplusn}{\sum}

\setlength{\textheight}{9.5in}

\newcommand{\Obs}{$\mathcal{O}$}
\newcommand{\OBS}{\mathcal{O}}

\title{Detection of a very serious error in the paper: "On identifiability of nonlinear ODE models and applications in viral dynamics"}

\author{Agostino Martinelli% <-this % stops a space
\thanks{A. Martinelli is with INRIA Rhone Alpes,
Montbonnot, France e-mail: {\tt agostino.martinelli@inria.fr}} }

\maketitle

\tableofcontents

\begin{abstract}
This erratum highlights
a very serious error in a paper published by SIAM Review in 2011.
The error is in Section 6.2 of  \cite{Miao11}. It is very important 
to notify this error because of the following two reasons: (i)
\cite{Miao11} is one of the most cited contributions in the field of identifiability of viral dynamics models, and (ii)
the error is relevant because, as a result of it, a very popular viral model (perhaps the most popular in the field of HIV dynamics) has been classified as 
identifiable. In contrast, three of its parameters are not identifiable, even locally.
This erratum first proves the non uniqueness of the three unidentifiable parameters by exhibiting infinitely many distinct but indistinguishable values of them. The non uniqueness is even local.
Then, this erratum details the error made by the authors of \cite{Miao11} which produced the claimed  (but false) local identifiability of all the model parameters.
\end{abstract}

\section{Introduction}\label{ChapterIntroduction}\label{ChapterIntroduction}
Section 6.2 of \cite{Miao11} provided an identifiability analysis of a very popular viral model that describes the HIV/AIDS dynamics. The model is characterized by six parameters ($\lambda$, $\delta$, $\eta$, $\rho$, $N$, $c$) and the analysis concerned the case when one of them ($\eta$) is time varying. 
The model is characterized by the following equations:

\begin{equation}\label{EquationHIVSystem}
\left\{\begin{array}{ll}
\dot{T}_U &= \lambda -\rho T_U -\eta(t) T_U V\\
\dot{T}_I &= \eta(t) T_U V -\delta T_I\\
\dot{V} &= N\delta T_I - cV\\
y &= [V,~T_U+T_I], \\
\end{array}\right.
\end{equation}
which provide the dynamics and the outputs. 
The above equations are Equations (6.15-6.17) in Section 6.2 of \cite{Miao11}. The reader is addressed to Section 6.2 of \cite{Miao11} for the definition of all the quantities that appear in these equations.

Unfortunately, the result of the identifiability analysis in Section 6.2 of \cite{Miao11} is wrong. Specifically, when $\eta$ is time varying, it is not true that all the model parameters are locally identifiable, as concluded by the authors of \cite{Miao11} at the end of Section 6.2. The result of a thorough identifiability analysis shows that only three of the five constant parameters are locally identifiable. The remaining two constant parameters ($\delta$, and $N$), and the time varying parameter ($\eta$) cannot be uniquely identified, even locally. 
In this paper, we prove the non uniqueness of $\delta$, $N$, and $\eta$ by exhibiting infinitely many distinct but indistinguishable values of these parameters (Section \ref{SectionNonUniqueness}). The non uniqueness is even local.
Finally,  in Section \ref{SectionError}
we detail the technical error in Section 6.2 of \cite{Miao11} that produced the claimed  (but false) local identifiability of all the model parameters.

\section{Non uniqueness of $\delta$, $N$, and $\eta$}\label{SectionNonUniqueness}

Let us refer to the viral model in Section  6.2 of \cite{Miao11}.
The constant parameters  $\lambda$, $\delta$, $\rho$, $N$, $c$, the three states $T_U(t), ~T_I(t),~V(t)$, and the time varying parameter $\eta(t)$ satisfy
the three dynamical equations in (\ref{EquationHIVSystem}).
Starting from the above parameters and states, we introduce the following new constant parameters, the following new states, and the following new time varying parameter, for any $\tau\in\mathbb{R}$:
\begin{equation}\label{EquationParametersTransformed}
\left\{\begin{array}{lll}
\lambda'& &= \lambda\\
%\delta'&=\delta'(\tau)&= \delta \rho\left/\left(( \rho- \delta)e^{ \rho\tau}+ \delta\right)\right.\\
\delta'&=\delta'(\tau)&= \frac{\delta \rho}{( \rho- \delta)e^{ \rho\tau}+ \delta}\\
\rho'& &= \rho\\
N'&=N'(\tau)&= Ne^{ \rho\tau}\\
c'& &= c, \\
\end{array}\right.
\end{equation}

\begin{equation}\label{EquationStatesTransformed}
\left\{\begin{array}{lll}
T_U'&=T_U'(t,~\tau)&=T_U(t)+ T_I(t)-\frac{ T_I(t)}{ \rho}\left(
 \delta e^{- \rho\tau}+ \rho- \delta
\right)\\
 T_I'&= T_I'(t,~\tau)&=\frac{ T_I(t)}{ \rho}\left(
 \delta e^{- \rho\tau}+ \rho- \delta
\right)\\
V'&=V'(t) &=V(t)\\
\end{array}\right.
\end{equation}

\begin{equation}\label{EquationEtaTransformed}
\begin{array}{ll}
\eta'&=\eta'(t,~\tau)=\frac{ \eta(t) T_U(t)  V(t)  \rho e^{ \rho \tau} +\left[T_I(t)  \delta^2-  T_I(t)  \delta  \rho-  \eta(t) T_U(t)  V(t)  \delta  \right] \left[e^{ \rho \tau}-1 \right]}{ V(t) \left[ T_I(t)\delta +  T_U(t)  \rho \right] e^{ \rho \tau} -  V(t)T_I(t)  \delta}\\
\end{array}
\end{equation}

We claim that all the above new constant parameters, the new states, and the new time varying parameter, are indistinguishable from the original ones, and this holds for any choice of $\tau\in\mathbb{R}$.

Our statement can be easily checked by a direct substitution.

First, the outputs $y_1(t)=T_U(t)+T_I(t)$ and $y_2(t)=V(t)$ remain the same because $T_U'(t,~\tau)+T_I'(t,~\tau)=T_U(t)+T_I(t)$ and $V'(t)=V(t)$.

Second, from (\ref{EquationHIVSystem}) and (\ref{EquationParametersTransformed}-\ref{EquationEtaTransformed}) it is possible to check, by a simple substitution, the validity of the following dynamics of the new three states (see all the details in Appendix \ref{Appendix}):

\begin{equation}\label{EquationDynamicsTransformed}
\left\{\begin{array}{ll}
\frac{d}{dt}T_U'(t,~\tau) &=\lambda'-\rho'T_U'(t,~\tau)-\eta'(t,~\tau)T_U'(t,~\tau)V'(t)\\
\frac{d}{dt}T_I'(t,~\tau) &=\eta'(t,~\tau)T_U'(t,~\tau)V'(t)-\delta'T_I'(t,~\tau)\\
\frac{d}{dt}V'(t) &=N'\delta'T_I'(t,~\tau)-c'V'(t)\\
\end{array}\right.
\end{equation}

In other words, the new three states, satisfy exactly the same dynamical equations in (\ref{EquationHIVSystem}) with the new parameters.

We conclude that
$T_U'$, $T_I'$, $V'$, 
$\eta'$, $\lambda'$, $\delta'$, $\rho'$, $N'$, and $c'$ cannot be distinguished from $T_U$, $T_I$, $V$, $\eta$, $\lambda$, $\delta$, $\rho$, $N$, and $c$, respectively. 
Regarding the model parameters, 
as $\eta'\neq\eta$, $\delta'\neq\delta$, and $N'\neq N$, we conclude that $\eta$, $\delta$, and $N$ are not identifiable. Finally, we remark that $\tau$ can take infinitesimal values and that, in the limit of $\tau\rightarrow0$, we have: 
$T_U'\rightarrow T_U$, $T_I'\rightarrow T_I$, $V'\rightarrow V$,
$\eta'\rightarrow\eta$, $\delta'\rightarrow\delta$, and $N'\rightarrow N$. Therefore, we conclude that the three parameters $\eta$, $\delta$, and $N$ are not even locally identifiable.

In Chapter 6 of \cite{ArXiv22}, we also provided a numerical test for the validity of the above fundamental new result (Section 6.6 of \cite{ArXiv22}). In particular, we adopted the same data set introduced in \cite{Villa19b}. Finally, in Section 6.7 of \cite{ArXiv22} we characterized the minimal information that must be added to the knowledge of the two outputs, in order to make identifiable all the model parameters (and to make observable all the states).

We wish to emphasize that our claim was here easily proved by a simple substitution. What is not trivial is the analytic derivation of the expressions in 
(\ref{EquationParametersTransformed}-\ref{EquationEtaTransformed}). These expressions were analytically determined in \cite{ArXiv22}, starting from the analytic solution
of a fundamental problem, strongly connected with the problem of determining the identifiability of time varying parameters. This is the unknown input observability problem, whose general analytic solution was recently introduced in \cite{IF22}.
The analytic solution given in \cite{ArXiv22} is strongly based on the solution in \cite{IF22} and provides the identifiability of any ODE model, even in the presence of time varying parameters and in the presence of any type of nonlinearity (and not necessarily polynomial, as for the methods discussed in \cite{Miao11}). 
Additionally, in contrast with the methods discussed in \cite{Miao11}, it works automatically (i.e., by following the steps of a systematic procedure).
Finally, in the presence of unidentifiability, it also provides, automatically, the local indistinguishable states and parameters (i.e., for the specific case, the non trivial expressions in 
(\ref{EquationParametersTransformed}-\ref{EquationEtaTransformed}), for infinitesimal $\tau$).

\section{The error in Section 6.2 of \cite{Miao11}}\label{SectionError}

The derivation given in Section 6.2 of \cite{Miao11} contains a very serious error that will be detailed in Section \ref{SubSectionError}. To explain this very serious error, we first need to highlight another minor error that could simply be a typo (Section \ref{SubSectionTypo}).

\subsection{Two typos in Equation (6.23) in \cite{Miao11}}\label{SubSectionTypo}

%First, Equation (6.23) in \cite{Miao11} contains two errors (probably typos). 
Equation (6.23) was obtained from (6.21) and (6.22), which are correct. From (6.21) and (6.22), we obtain:

\begin{equation}\label{Equation6.23Corrected}
\ddot{y}_1 y_2 \dot{y}_2-\dot{y}_1 y_2 \ddot{y}_2-\delta y_1 y_2 \ddot{y}_2+\lambda y_2 \ddot{y}_2-(\delta+c) \dot{y}_1 y_2 \dot{y}_2+(\delta \rho-\delta^2-\delta c) y_1 y_2 \dot{y}_2
\end{equation}
\[
+(\rho+\delta) \dot{y}_1 y_2 \dot{y}_2+\lambda c y_2 \dot{y}_2+\rho c \dot{y}_1 y_2^2+(\rho \delta c-\delta^2 c) y_1 y_2^2-N \delta y_1 \ddot{y}_1 y_2
\]
\[
+c \ddot{y}_1 y_2^2-N \delta (\rho+\delta) y_1 \dot{y}_1 y_2-N \delta^2 \rho y_1^2 y_2+N \delta^2 \lambda y_1 y_2=0
\]

The differences between (\ref{Equation6.23Corrected}) and (6.23) are the following two:

\begin{enumerate}

\item The sixth term in (6.23) is not $(\delta \rho+\rho+\delta-\delta^2-\delta c) y_1 y_2 \dot{y}_2$ but $(\delta \rho-\delta^2-\delta c) y_1 y_2 \dot{y}_2+(\rho+\delta) \dot{y}_1 y_2 \dot{y}_2$.

\item The seventh  term in (6.23) is not $c y_2 \dot{y}_2$ but $\lambda c y_2 \dot{y}_2$.

\end{enumerate}

\subsection{The very serious error made by the authors of \cite{Miao11}}\label{SubSectionError}

Let us denote the expression in (\ref{Equation6.23Corrected}) by $\phi=\phi(y_1,\dot{y}_1,\ddot{y}_1,y_2,\dot{y}_2,\ddot{y}_2,\lambda,\delta,\rho,c,N)$. 
In this notation, Equation (\ref{Equation6.23Corrected}) becomes:
\[
\phi=0.
\]

In Section 6.2 of  \cite{Miao11}, the local identifiability of all the parameters was claimed by referring to Equation (6.23). In particular, the authors suggested two methods to obtain the local identifiability: (i) by using the third statement of Theorem 3.8 in \cite{Miao11}, and (ii) by using the implicit function Theorem (in accordance with the procedure given in Section 3.4 of \cite{Miao11}). We cannot comment the first method because the authors only provided the last polynomial 
(i.e., $B_n(u,y,\theta_1,\ldots,\theta_q)$) of Theorem 3.8. They did not provide the remaining  polynomials mentioned by the theorem (i.e., $B_1,\ldots,B_{n-1}$).

Regarding the second method, the error is pointed out as follows.
In accordance with the procedure in Section 3.4 of \cite{Miao11}, we need, first of all, to compute the first 4 time derivatives of $\phi$. In this manner, we obtain four additional equations, i.e.,
\[
\phi^1\equiv\frac{d}{dt}\phi=0,~~\ldots,~~\phi^4\equiv\frac{d^4}{dt^4}\phi=0
\]
Then, we include all the above functions in the vector function $\Phi$. In other words, we introduce the vector function $\Phi\equiv\left[
 \phi,~ \phi^1, ~\phi^2, ~\phi^3, ~\phi^4 \right]^T$. As the parameters $\lambda,\delta,\rho,c,N$ are constant, $\Phi$ is a function of them and of the two outputs ($y_1,~y_2$) and their time derivatives up to the sixth order. 
To check the identifiability of $\lambda,\delta,\rho,c,N$, the procedure in Section 3.4 of \cite{Miao11} uses the following equation:
 
 \begin{equation}\label{EquationVectorForm}
 \Phi=0
 \end{equation}

Now, by using the implicit function Theorem, we need to compute the Jacobian of $\Phi$ with respect to $\lambda,\delta,\rho,c,N$. By proceeding in this manner, we obtain a $5\times5$ matrix that we denote by $M$. The entries of $M$ are expressed in terms of $\lambda,\delta,\rho,c,N,y_1,y_2$, and the time derivatives of $y_1,y_2$, up to the sixth order. 

{\bf The very serious error made by the authors of \cite{Miao11}} was to consider INDEPENDENT all the above quantities, i.e., $\lambda,\delta,\rho,c,N,y_1,y_2$, and the time derivatives of $y_1,y_2$, up to the sixth order. By doing this we obtain that the rank of $M$ is equal to 5 (and the authors concluded that all the model parameters are identifiable).

On the other hand,
the two output functions, $y_1(t), ~y_2(t)$, and their time derivatives up to any order, cannot take arbitrary values . {\bf They must satisfy the constraints due to the system dynamics}.
%In \cite{Miao11}, neither in Section 3.4, nor in Section 6.2, this fundamental aspect was mentioned.
In accordance with the dynamics constraints, we must express the two outputs ($y_1,~y_2$) and their time derivatives up to the sixth order, in terms of the state ($T_U, ~T_I, ~V$), the model parameters $\lambda,\delta,\rho,c,N$, and the time varying parameter $\eta$ and its time derivatives up to the fifth order. We obtain all these expressions by repetitively using (\ref{EquationHIVSystem}). Finally, we compute the rank of $M$ and, by a direct computation, we obtain the following fundamental property. Independently of the values of $\lambda,\delta,\rho,c,N,T_U,T_I,V,\eta,\eta^{(1)},\eta^{(2)},\eta^{(3)},\eta^{(4)},\eta^{(5)}$ (with $\eta^{(k)}\equiv\frac{d^k\eta}{dt^k}$), this rank cannot exceed $4<5$. As a result, it is very erroneous to conclude that Equation (\ref{EquationVectorForm}) allows us to locally identify all the five parameters $\lambda,\delta,\rho,c,N$.
%We remark that, if $\phi$ were defined by the erroneous expression in (6.23) (instead of the correct expression in (\ref{Equation6.23Corrected})), the aforementioned rank would become $5$.

\section{Conclusion}

This erratum highlighted a very serious error on a paper
published by SIAM Review in 2011.
The error is in Section 6.2 of  \cite{Miao11}. 

The error is relevant because, as a result of it, a very popular viral model (perhaps the most popular in the field of HIV dynamics) has been classified as 
identifiable. In contrast, three of its parameters are not identifiable, even locally.

This erratum first proved the non uniqueness of the three unidentifiable parameters by exhibiting infinitely many  distinct but indistinguishable values of them (Section \ref{SectionNonUniqueness}). The non uniqueness is even local.
Then, this erratum detailed the very serious error made by the authors of \cite{Miao11} which produced the claimed  (but false) local identifiability of all the model parameters (Section \ref{SubSectionError}).

\newpage

\appendix

\section{Technical details for the derivation of Eq. (\ref{EquationDynamicsTransformed})}\label{Appendix}

\subsection*{First equation}
Let us prove the validity of the first equation in (\ref{EquationDynamicsTransformed}), i.e., the validity of
\[
\frac{d}{dt}T_U' =\lambda'-\rho'T_U'-\eta'T_U'V'
\]
We prove the equality of the two members.

\subsubsection*{First member}
From the first equation in (\ref{EquationStatesTransformed}), we obtain:
\[
\frac{d}{dt}T_U' = \frac{d}{dt}T_U+ \frac{d}{dt}T_I-\frac{ \frac{d}{dt}T_I}{ \rho}\left(\delta e^{- \rho\tau}+ \rho- \delta
\right).
\]

By using the first two equations in (\ref{EquationHIVSystem}), we obtain:

\begin{equation}\label{EquationFirstEqFirstM}
\frac{d}{dt}T_U' = \lambda -\rho T_U -\eta T_U V+  \eta T_U V -\delta T_I-\frac{  \eta T_U V -\delta T_I}{ \rho}\left(\delta e^{- \rho\tau}+ \rho- \delta
\right)=
\end{equation}
\[
=\lambda -\rho T_U  -\delta T_I+\frac{  \delta T_I-\eta T_U V}{ \rho}\left(\delta e^{- \rho\tau}+ \rho- \delta\right).
\]

\subsubsection*{Second member}

The second member is $\lambda'-\rho'T_U'-\eta'T_U'V'$. By using (\ref{EquationParametersTransformed}), we obtain:

\[
\lambda'-\rho'T_U'-\eta'T_U'V'=\lambda-\rho T_U'-\eta'T_U'V'
\]

By using the expression of $T_U'$ and $V'$ in (\ref{EquationStatesTransformed}), and by using the expression of $\eta'$ in (\ref{EquationEtaTransformed}), after some simplifications, we obtain:

\begin{equation}\label{EquationFirstEqSecondM}
\lambda'-\rho'T_U'-\eta'T_U'V'=
-\frac{T_I \delta^2 - \lambda  \rho + T_U \rho^2 - T_I \delta^2 e^{-\rho \tau} - T_U V \delta  \eta  + T_U V \eta  \rho + T_U V \delta  \eta  e^{-\rho \tau}}{\rho}
\end{equation}

\vskip.1cm
By a direct computation it is easy to check that the two members (given in (\ref{EquationFirstEqFirstM}) and (\ref{EquationFirstEqSecondM}), respectively) coincide, i.e., the following equality holds

\[
\lambda -\rho T_U  -\delta T_I+\frac{  \delta T_I-\eta T_U V}{ \rho}\left(\delta e^{- \rho\tau}+ \rho- \delta\right)
=
-\frac{T_I \delta^2 - \lambda  \rho + T_U \rho^2 - T_I \delta^2 e^{-\rho \tau} - T_U V \delta  \eta  + T_U V \eta  \rho + T_U V \delta  \eta  e^{-\rho \tau}}{\rho}
\]

%\newpage
%
%\[
%\lambda -\rho T_U  -\delta T_I+\frac{  \delta T_I-\eta T_U V}{ \rho}\left(\delta e^{- \rho\tau}+ \rho- \delta\right)
%+\frac{T_I \delta^2 - \lambda  \rho + T_U \rho^2 - T_I \delta^2 e^{-\rho \tau} - T_U V \delta  \eta  + T_U V \eta  \rho + T_U V \delta  \eta  e^{-\rho \tau}}{\rho}
%\]
%
%\[
%\lambda\rho -\rho^2 T_U  -\delta\rho T_I+(\delta T_I-\eta T_U V)\left(\delta e^{- \rho\tau}+ \rho- \delta\right)
%+
%\]
%\[
%T_I \delta^2 - \lambda  \rho + T_U \rho^2 - T_I \delta^2 e^{-\rho \tau} - T_U V \delta  \eta  + T_U V \eta  \rho + T_U V \delta  \eta  e^{-\rho \tau}
%\]
%
%\[
% -\rho^2 T_U  -\delta\rho T_I+(\delta T_I-\eta T_U V)\left( \rho- \delta\right)
%+
%T_I \delta^2  + T_U \rho^2  - T_U V \delta  \eta  + T_U V \eta  \rho 
%\]
%
%
%\[
%  -\delta\rho T_I+(\delta T_I-\eta T_U V)\left( \rho- \delta\right)
%+
%T_I \delta^2    - T_U V \delta  \eta  + T_U V \eta  \rho 
%\]
%
%\[
% (-\eta T_U V)\left( \rho- \delta\right)
% - T_U V \delta  \eta  + T_U V \eta  \rho 
%\]
%
%\[
% -\left( \rho- \delta\right)
% -  \delta    +   \rho =0
%\]
%
%\newpage
%

\subsection*{Second equation}
Let us prove the validity of the second equation in (\ref{EquationDynamicsTransformed}), i.e., the validity of
\[
\frac{d}{dt}T_I' =\eta'T_U'V'-\delta'T_I'
\]

We prove the equality of the two members.

\subsubsection*{First member}
From the second equation in (\ref{EquationStatesTransformed}), we obtain:
\[
\frac{d}{dt}T_I'=
\frac{\frac{d}{dt} T_I}{ \rho}\left(
 \delta e^{- \rho\tau}+ \rho- \delta
\right).
\]

From the second equation in (\ref{EquationHIVSystem}), we obtain:

\begin{equation}\label{EquationSecondEqFirstM}
\frac{d}{dt}T_I'=
\frac{ \eta T_U V -\delta T_I}{ \rho}\left(
 \delta e^{- \rho\tau}+ \rho- \delta
\right).
\end{equation}

\subsubsection*{Second member}

The second member is $\eta'T_U'V'-\delta'T_I'$. By using (\ref{EquationParametersTransformed}), we obtain:

\[
\eta'T_U'V'-\delta'T_I'=
\eta'T_U'V'-\frac{\delta \rho}{( \rho- \delta)e^{ \rho\tau}+ \delta}T_I'
\]

By using the expression of $T_U'$, $T_I'$, and $V'$ in (\ref{EquationStatesTransformed}), and by using the expression of $\eta'$ in (\ref{EquationEtaTransformed}), after some simplifications, we obtain:

\begin{equation}\label{EquationSecondEqSecondM}
\eta'T_U'V'-\delta'T_I'=
-\frac{e^{-\rho\tau} (T_I \delta - T_U V \eta) (\delta - \delta  e^{\rho\tau} + \rho  e^{\rho\tau})}{\rho}
\end{equation}

\vskip.1cm
By a direct computation it is easy to check that the two members (given in (\ref{EquationSecondEqFirstM}) and (\ref{EquationSecondEqSecondM}), respectively) coincide, i.e., the following equality holds

\[
\frac{ \eta T_U V -\delta T_I}{ \rho}\left(
 \delta e^{- \rho\tau}+ \rho- \delta
\right)
=
-\frac{e^{-\rho\tau} (T_I \delta - T_U V \eta) (\delta - \delta  e^{\rho\tau} + \rho  e^{\rho\tau})}{\rho}
\]

%
%\newpage
%
%
%\[
%(\eta T_U V -\delta T_I)\left(
% \delta e^{- \rho\tau}+ \rho- \delta
%\right)
%+ (T_I \delta - T_U V \eta) (\delta e^{-\rho\tau} - \delta   + \rho  )
%\]
%
%
%
%\newpage
%
%

\subsection*{Third equation}
Let us prove the validity of the third equation in (\ref{EquationDynamicsTransformed}), i.e., the validity of
\[
\frac{d}{dt}V' =N'\delta'T_I'-c'V'
\]

We prove the equality of the two members.

\subsubsection*{First member}
From the third equation in (\ref{EquationStatesTransformed}), we obtain:
\[
\frac{d}{dt}V'=
\frac{d}{dt}V
\]

From the third equation in (\ref{EquationHIVSystem}), we obtain:

\begin{equation}\label{EquationThirdEqFirstM}
\frac{d}{dt}V'=
 N\delta T_I - cV
\end{equation}

\subsubsection*{Second member}

The second member is $N'\delta'T_I'-c'V'$. By using (\ref{EquationParametersTransformed}), we obtain:

\[
N'\delta'T_I'-c'V'=
 Ne^{ \rho\tau}\frac{\delta \rho}{( \rho- \delta)e^{ \rho\tau}+ \delta}T_I'-cV'
\]

By using  (\ref{EquationStatesTransformed}), we obtain:

\begin{equation}\label{EquationThirdEqSecondM}
N'\delta'T_I'-c'V'=
 Ne^{ \rho\tau}\frac{\delta \rho}{( \rho- \delta)e^{ \rho\tau}+ \delta}
 \frac{ T_I}{ \rho}\left(
 \delta e^{- \rho\tau}+ \rho- \delta
\right)-cV
\end{equation}

\vskip.1cm
By a direct computation it is easy to check that the two members (given in (\ref{EquationThirdEqFirstM}) and (\ref{EquationThirdEqSecondM}), respectively) coincide, i.e., the following equality holds

\[
 N\delta T_I - cV=
Ne^{ \rho\tau}\frac{\delta \rho}{( \rho- \delta)e^{ \rho\tau}+ \delta}\frac{ T_I}{ \rho}\left(
 \delta e^{- \rho\tau}+ \rho- \delta
\right)
-cV
\]

%\newpage
%
%\[
% N\delta T_I 
%-Ne^{ \rho\tau}\frac{\delta \rho}{( \rho- \delta)e^{ \rho\tau}+ \delta}\frac{ T_I}{ \rho}\left(
% \delta e^{- \rho\tau}+ \rho- \delta
%\right)
%\]
%
%\[
%\rho N\delta T_I 
%-Ne^{ \rho\tau}\frac{\delta \rho}{( \rho- \delta)e^{ \rho\tau}+ \delta}T_I\left(
% \delta e^{- \rho\tau}+ \rho- \delta
%\right)
%\]
%
%
%\[
% \left(( \rho- \delta)e^{ \rho\tau}+ \delta\right)
%\rho N\delta T_I 
%-Ne^{ \rho\tau}
%\delta \rho T_I\left(
% \delta e^{- \rho\tau}+ \rho- \delta
%\right)
%\]
%
%\[
% \left(( \rho- \delta)e^{ \rho\tau}+ \delta\right)
%-
%e^{ \rho\tau}
%\left(
% \delta e^{- \rho\tau}+ \rho- \delta
%\right)
%\]
%
%\[
% \left(( \rho- \delta)e^{ \rho\tau}+ \delta\right)
%-
%\left(
% \delta + \rho e^{ \rho\tau}- \delta e^{ \rho\tau}
%\right)
%\]
%
%\newpage
%

%%%%%%%%%%%%%%%%%%%%%%%%%%%%%%%%%%%%%%%%%%%%%%%%
%%%%%%%%%%%%%%%%%%%%%%%%%%%%%%%%%%%%%%%%%%%%%%%%
%%%%%%%%%%%%%%%%%%%%%%%%%%%%%%%%%%%%%%%%%%%%%%%%
%%%%%%%%%%%%%%%%%%%%%%%%%%%%%%%%%%%%%%%%%%%%%%%%
%%%%%%%%%%%%%%%%%%%%%%%%%%%%%%%%%%%%%%%%%%%%%%%%

\end{document}